\documentclass[11pt]{article}
\usepackage{latexsym,amsmath,amsfonts,mathrsfs}
\newcommand{\nref}[1]{(\ref{#1})}
\def\Cal{\cal}

\usepackage[cp1251]{inputenc} 
\def\text#1{\hbox{#1}}

\def\endproof{\mbox{\ $\Box$}}

\newcommand{\R}{\mathbb R}
\newcommand{\N}{\mathbb N}

\newcommand{\Z}{\mathbb Z}

\def\1{\mbox{1\hspace{-.20em}I}}

\newcommand{\ud}{\mathrm{d}}

\newtheorem{theorem}{Theorem}[section]
\newtheorem{proposition}{Proposition}[section]

\newtheorem{remark}{Remark}[section]

\newcommand{\CY}{{\Cal{Y}}}

\newcommand{\CN}{{\Cal{N}}}

\def\a{\alpha}

\def\e{\varepsilon}
\def\b{\beta}
\def\g{\gamma}

\def\d_1{\gamma_1}
\def\l{\left}
\def\r{\right}

\def\Var{\rm {Var}}

\def\Var{{\rm Var}}

\voffset=-20mm \numberwithin{equation}{section}
\begin{document}

\title{Minimax nonparametric testing in a problem related to the Radon transform}

\author{Yuri I. Ingster\footnote{Research was partially supported
by RFBR grant 11-01-00577 and by the grant NSh--4472.2010.1}
\\ {\small Department of Mathematics II,} \\
{\small Saint Petersburg State Electrotechnical University,}\\
{\small  Russia}\\  \\ Theofanis Sapatinas\\ {\small Department of
Mathematics and
Statistics,}\\
{\small  University of Cyprus,}\\
{\small  Cyprus}\\ \\
 and\\ \\
 Irina A. Suslina
 \\ {\small Department of Mathematics,} \\
{\small Saint Petersburg State  University of Information Technologies,
Mechanics and Optics,}\\
{\small  Russia}}


\maketitle

\begin{abstract}
We consider the detection problem of a two-dimensional function from noisy observations of its integrals over lines.
We study both rate and sharp asymptotics for the error probabilities in the minimax setup.
By construction, the derived tests are non-adaptive. We also construct a minimax rate-optimal adaptive test of
rather simple structure.



\end{abstract}

\section{Introduction}
The problem of tomography is to reconstruct a two-dimensional function (image) from its Radon transform, i.e.,
from observations of its integrals over lines. This problem, and its extension to higher dimensions,
appears in different scientific fields such as radio astronomy and medical imaging
(see, e.g.,  \cite{Deans}, \cite{Her}, \cite{Nat}).  We consider the tomography problem from a
statistical perspective that can be formulated as a problem of reconstructing a two-dimensional
function from its noisy Radon transform (see, e.g., \cite{cav}, \cite{Cav-Tsyb}, \cite{KKLPP}, \cite{KLP}).

Despite some work on the minimax estimation problem of a two-dimensional function from its noisy Radon transform
(see \cite{Cav-Tsyb}, \cite{KLP}, \cite{KorTsy}), to the best of our knowledge, there exist no work on the
corresponding minimax detection problem. The general statement of this problem is given in Section 2,
while some preliminaries and notation in the minimax signal detection framework are presented in Section 3.
Within this framework, in Section 4, we consider the detection problem of a two-dimensional function from
its noisy Radon transform and study both rate and sharp asymptotics for the error probabilities.
By construction, the derived tests are non-adaptive. A rate-optimal adaptive test of rather simple
structure is also constructed. The proofs are given in the Appendix.

\section{Formulation of the problem}

\subsection{The Radon transform}
Denote by $\|\cdot\|$ the standard Euclidean norm in $\R^2$, i.e., $\|x\|=(x_1^2+x_2^2)^{1/2}$,
$x=(x_1,x_2) \in \R^2$. Let $H=\{x\in\R^2: \|x\| \leq 1\}$,  be the unit disk in $\R^2$, and
let $\mu$ denote the Lebesgue measure in $\R^2$.
Consider the integrals of a function $f:H \mapsto \R$ over all lines
that intersect $H$. The lines are parameterized by the length
$u\in[0,1]$ of the perpendicular from the origin to the line and by
the orientation $\varphi \in [0,2\pi)$ of this perpendicular.
Suppose that $f \in L^1(H,\mu) \cap
L^2(H,\mu)$. Define the Radon transform of the function $f$ by
\begin{equation}
\label{RDTr}
{\cal R}f(u,\varphi) = \frac{\pi}{2 \sqrt{1-u^2}} \int_{-\sqrt{1-u^2}}^{\sqrt{1-u^2}}
\,\,f(u \cos\varphi-t\sin\varphi, u \sin\varphi+t\cos\varphi ) \,\ud t, \quad (u,\varphi) \in S,
\end{equation}
where
$$
S=\{(u,\varphi): u\in [0,1], \, \varphi \in [0, 2\pi)\}.
$$
Thus, the Radon transform ${\cal R}f$ is $\pi$ times the average of $f$ over the
line segment (parametrized by $(u,\varphi)$) that intersects $H$. It is natural to
consider ${\cal R}f$ as an element of $L^2(S,\mu_0)$, where $\mu_0$ is the measure on $S$ defined
 by
$$
\ud \mu_0(u,\varphi) = \frac{2\, \sqrt{1-u^2}}{\pi}
\,\ud \varphi,\quad (u,\varphi) \in S.
$$

\subsection{The Gaussian white noise model}
Consider now the Gaussian white noise model
\begin{equation}
\label{gwmRDTFan}
\ud Y_{\e}(u,\varphi)={\cal R}f(u,\varphi)\, \ud u\,\ud \varphi+\e\, \ud W(u,\varphi),
 \quad (u,\varphi) \in S,
\end{equation}
where $W$ is a standard Wiener sheet on $S$ (i.e., the primitive of white noise on $S$)
 and $\e>0$ is a small parameter (the noise level).  Although this model
is continuous and real data are typically discretely sampled,  versions of it have been extensively studied
in the nonparametric literature and are considered as idealized models that provide, subject
to some limitations, approximations to many, sampled-data, nonparametric models
 (see, e.g., \cite{br-lo}, \cite{Candes-Donoho2002}, \cite{efr-sam}, \cite{nus}).

The Gaussian white noise model (\ref{gwmRDTFan}) may also seem initially rather remote.
One may, however, be
helped by the observation that what it really means is the following: for any function
$g \in L^2(S,\mu_0)$, the integral
$$
\iint_{S} g(u,\varphi){\cal R}f(u,\varphi)\, \ud u \,\ud \varphi
$$
can be observed with Gaussian error having zero mean and variance equal to
$\e^2 \iint_{S}g^2(u,\varphi)\,\ud u\,\ud \varphi$ (see, e.g., \cite{Candes-Donoho2002}).

The Radon transform  ${\cal R}$ is a compact operator and its
singular value decomposition (SVD)
 is well-known
(see, e.g., \cite{Nat}). To introduce it, let $\N=\{1,2,\ldots\}$ be
the set of natural numbers, set  $\Z_+=\N\cup\{0\}$, and define a set of double indices giving
rise to the following lattice quadrant
\begin{equation}
\label{GammaFF}
\Gamma=\{\nu: \nu=(j,l),\,\, j,l \in \Z_+\}.
\end{equation}
An orthonormal complex-valued basis for $L^2(H,\mu)$ is given by
\begin{equation}
\label{complF1} \tilde{\phi}_{\nu}(r,\theta)=\pi^{-1/2}
(j+l+1)^{1/2} Z_{j+l}^{|j-l|}(r) \exp\{i(j-l)\theta\}, \quad \nu\in
\Gamma,
\end{equation}
where $x=(r\cos\theta, r\sin\theta) \in H,$ with $Z_a^b$
denoting the Zernike polynomial of degree $a$ and order $b$, with
$a, b \in \Z_+$ (see, e.g., \cite{Deans}). The corresponding
orthonormal complex-valued basis in $L^2(S,\mu_0)$ is
\begin{equation}
\label{complF2} \tilde{\psi}_{\nu}(u,\varphi)=\pi^{-1/2} U_{j+l}(u)
\exp\{i(j-l)\varphi\}, \quad \nu\in \Gamma, \quad (u,\varphi) \in S,
\end{equation}
where
$$
U_m(\cos\theta)=\frac{\sin((m+1)\theta)}{\sin\theta}, \quad m \in \Z_+,
\quad \theta \in [0, 2\pi),
$$
are the Chebyshev polynomials of the second kind. We then have (see, e.g., \cite{Cav-Tsyb})
$$
{\cal R}\tilde{\phi}_{\nu} = b_{\nu}\tilde{\psi}_{\nu}
$$
with singular values
\begin{equation}
\label{radsing} b_{\nu}=\pi(j+l+1)^{-1/2}, \quad \nu \in \Gamma.
\end{equation}

Since we work with real-valued functions $f$, the complex-valued bases (\ref{complF1}) and
(\ref{complF2}) are identified, in standard fashion, with the equivalent  real-valued
orthonormal bases $\phi_{\nu}$ and $\psi_{\nu}$, $\nu \in \Gamma$, respectively, defined by
\begin{equation}
\label{Fourorth}
\phi_{\nu}=
\begin{cases}
&\sqrt{2}\, \text{Re}(\tilde{\phi}_{\nu}), \quad \text{if} \ j>l,\cr
&\tilde{\phi}_{\nu}, \quad
\quad \quad \quad \,\,\,
\text{if} \ j=l,\cr
&\sqrt{2}\,  \text{Im}(\tilde{\phi}_{\nu}), \quad \text{if} \ j<l,
\end{cases}
\end{equation}
with an analogous expression for $\psi_{\nu}$, $\nu \in \Gamma$.

Hence,  by standard calculations (see, e.g., \cite{Cav-Tsyb}, \cite{ISS}) and an application of the spectral theorem
for the self-adjoint compact operator ${\cal R}^* {\cal R}$ (${\cal R}^*$ being the adjoint of
 ${\cal R}$), the Gaussian white noise model (\ref{gwmRDTFan})  generates the following equivalent discrete
 observational model in the Fourier
 domain, called the Gaussian sequence model,
\begin{equation}
\label{gsmRDtr}
y_{\nu}=b_{\nu} \theta_{\nu}+\e\, \xi_{\nu}, \quad \nu \in \Gamma,
\end{equation}
where $y_{\nu}=\langle {\cal R}f, \psi_{\nu}\rangle$, $\nu \in \Gamma$, are the ``observations'',
$b_{\nu}$, $\nu \in \Gamma$, are the singular values of the Radon operator ${\cal R}$, given by (\ref{radsing}),
$\theta_{\nu}=\langle f, \phi_{\nu} \rangle$, $\nu \in \Gamma$, are the Fourier coefficients of $f$ with respect
to $\phi_{\nu}$, given by (\ref{Fourorth}), and $\xi_{\nu}$, $\nu\in \Gamma$, are independent and identically
distributed (iid) standard Gaussian random variables, i.e., $\xi_{\nu} \stackrel{iid}{\sim} \CN(0,1)$, $\nu \in\Gamma$.

\subsection{The class of functions}
Crucial to the suggested detection methodology is the idea of considering minimax detection over certain classes
of functions in  $f \in L^2(H,\mu)$. Following \cite{Cav-Tsyb}, we consider a special class of functions with
polynomially decreasing coefficients $\theta=\{\theta_{\nu}\}_{\nu \in \Gamma}$, i.e., for some $p>0$, $L>0$,
\begin{equation}
\label{functClF}
\mathcal{F}(p,L)=\left\{f=\sum_{\nu\in\Gamma}\theta_{\nu}\phi_{\nu}:\
\ \theta\in\tilde\Theta(p,L)\right\}
\end{equation}
with
\begin{equation}
\label{functthlF}
\tilde\Theta(p,L)=\left\{\theta \in l^2:\ \ \sum_{\nu \in \Gamma,\, \nu \neq (0,0)}
(j+1)^{2 p}(l+1)^{2 p}\theta_{\nu}^2\leq L^2\right\}.
\end{equation}

It has been shown that $\mathcal{F}(p,L)$ can be identified with the set of functions $f$ which have $2p$ weak
derivatives (provided $2p$ is an integer) that are squared-integrable on $H$ with respect to the modified
dominating measure
$$
\ud \mu_{2 p +1}(x)=(1-\|x\|^2)^{2 p}\,\ud \mu(x), \quad x \in H.
$$
This is weaker than the square-integrability with respect to $\mu$ assumed for the usual Sobolev spaces
(see Proposition 2.2 in \cite{John-Silv}).

\subsection{The aim}

The goal is to determine whether the two-dimensional function $f$
corresponds to a known ``etalon" function $f_0$ (i.e., to test the null hypothesis $H_0: f=f_0$) or there exists a
difference between $f$ and $f_0$ (i.e., against the alternative hypothesis $H_1: f=f_0+\Delta f$ with $\Delta f \in \mathcal{F}(p,L)$
(see (\ref{functClF})--(\ref{functthlF})), based on the observation of a trajectory
$\{Y_{\e}=Y_{\e}(u,\varphi)\}$, $(u,\varphi) \in S$, from the
Guassian white noise model (\ref{gwmRDTFan}).

From mathematical point of view, we can take $f_0=0$ by passing to
the observation $\tilde Y_\e$ with $\ud\tilde Y_\e=\ud Y_\e(u,\phi)
-{\cal R}f_0\,\ud u\,\ud \phi$. For this reason, without loss of
generality,  we assume in the sequel that $f_0=0$, use $f$ in place
of $\Delta f$, and take the observation $Y_\e$. In order to avoid
having a trivial power (see below), our ultimate goal is to
determine whether $f$ satisfies (\ref{1.3}) (see below), using only
tests calibrated in such a way that if one had run them in the
absence of an $f \in \mathcal{F}(p,L)$, a certain restriction of the
significance level (error probability) is met.

In the sequel, we elaborate on the set under the alternative hypothesis and the suggested test statistics
that provide a good quality of testing in the minimax framework. Before going into the details, however,
we give the necessary preliminaries on the minimax signal detection framework in the standard Gaussian
white noise model which provide the avenue for developing the suggested detection methodology and
deriving theoretical results for detecting a two-dimensional function from its noisy Radon transform.

Hereafter, the relation $A_{\e} \sim B_{\e}$ means that $A_{\e}/B_{\e} \to 1$
as $\e \to 0$  while the relation $A_{\e} \asymp B_{\e}$ means that there
exists absolute constants $0 <c_1 \leq c_2 <\infty$ and $\e_0>0$ small
enough such that $c_1 \leq A_{\e}/B_{\e} \leq c_2$ for $0<\e \leq  \e_0$.

\section{Signal detection in the Gaussian sequence model: the minimax framework}\label{minimax}

Consider the Gaussian sequence model \nref{gsmRDtr}.  In order to avoid having a trivial minimax
hypothesis testing problem (i.e., trivial power), one usually needs to remove a neighborhood around
the functional parameter under the null hypothesis and to add some additional constraints, that are
typically expressed in the form of some regularity conditions, such as constraints on the derivatives,
of the unknown functional parameter of interest (see, e.g., \cite{IS.02}, Sections 1.3-1.4).

In view of the above observation, the main object of our study is the hypothesis testing problem
\begin{equation}\label{1.3}
H_0: \ \ \theta=0\quad\text{versus}\quad H_1:\ \
\sum_{\nu\in\Gamma}a_{\nu}^2\theta_{\nu}^2\leq
1,\quad\sum_{\nu\in\Gamma}\theta_{\nu}^2\ge r_\varepsilon^2,
\end{equation}
where  $\theta=\{\theta_{\nu}\}_{{\nu}\in\Gamma} \in l^2$,
$a_{\nu}\geq 0$, $\nu\in \Gamma$, and $r_\e>0,\ r_\e\to 0$,
 is a given family. It means that the
set under the alternative corresponds to an ellipsoid of semi-axes
$1/a_{\nu}$, $\nu\in\Gamma$, with an  $l^2$-ball of radius $r_\e$ removed.
(Here, $l^2=\{\zeta:\, \sum_{\nu\in\Gamma} \zeta^2_{\nu}<\infty\}$ with $\Gamma$ given by \nref{GammaFF}.)

Consider now the sequence $\eta=\{\eta_\nu\}_{\nu\in\Gamma}$ with
elements $\eta_{\nu}=\theta_{\nu}/\sigma_{\nu}$ where we set
$\sigma_{\nu}=1/b_{\nu}$, $\nu\in\Gamma$. In view of \nref{radsing},
the sequence $\eta=\{\eta_{\nu}\}_{\nu\in\Gamma} \in l_2$ and the
Gaussian sequence model \nref{gsmRDtr} takes the form
\begin{equation}
\label{ffgsm} y_{\nu}=\eta_{\nu}+\e \xi_{\nu}, \quad {\nu}\in\Gamma.
\end{equation}
The hypothesis testing problem
\nref{1.3} can now be written in the following equivalent form
\begin{equation}\label{1.4}
H_0: \ \ \eta=0\quad\text{versus}\quad H_1:\ \eta\in \Theta(r_\e),
\end{equation}
where the set under the alternative, i.e., $\Theta(r_\e)$, is
determined by the constraints
\begin{equation}\label{ell}
\Theta=\{\eta\in l^2: \ \sum_{\nu\in\Gamma}a_{\nu}^2\sigma_{\nu}^2\eta_{\nu}^2\leq
1\},\quad\Theta(r_\e)=\{\eta\in\Theta:\
\sum_{\nu\in\Gamma}\sigma_{\nu}^2\eta_{\nu}^2\ge r_\varepsilon^2\},
\end{equation}
i.e., the set under the alternative corresponds to an ellipsoid of
semi-axes $1/(a_{\nu}\sigma_{\nu})$, $\nu\in\Gamma$, with an
ellipsoid of semi-axes $r_\varepsilon/\sigma_{\nu}$, $\nu\in\Gamma$,
removed.

We are therefore interesting in the minimax efficiency of the
hypothesis testing problem \nref{1.4}-\nref{ell} for a given family
of sets $\Theta_{\e}=\Theta(r_\e)  \subset l^2$. It is characterized
by asymptotics, as $\e\to 0$, of the minimax error probabilities in
the problem at hand. Namely, for a (randomized) test $\psi$ (i.e., a
measurable function of the observation $y=\{y_\nu\}_{\nu\in\Gamma}$
taking values in $[0,1]$), the null hypothesis is rejected with
probability $\psi(y)$ and is accepted with probability $1-\psi(y)$.
Let $P_{\e,\eta}$ be the probability measure for the Gaussian
sequence model \nref{ffgsm} and denote by $E_{\e,\eta}$ the
expectation over this probability measure. Let $\a_{\e}(\psi) =
E_{\e,0}\psi$ be its type I error probability, and let $
\b_{\e}(\Theta_{\e},\psi) = \sup_{\eta\in
\Theta_{\e}}E_{\e,\eta}(1-\psi) $ be its maximal type II error
probability. We consider two criteria of asymptotic optimality:

(1) The first one corresponds to the classical Neyman-Pearson
criterion. For $\a\in (0,1)$, we set
$$
\b_{\e}(\Theta_{\e},\a) = \inf_{\psi:\,
\a_{\e}(\psi)\le\a}\b_{\e}(\Theta_{\e},\psi).
$$
We call a family of tests $\psi_{\e,\a}$ {\it asymptotically
minimax} if
$$
\a_{\e}(\psi_{\e,\a})\le\a+o(1),\quad
\b_{\e}(\Theta_{\e},\psi_{\e,\a}) = \b_{\e}(\Theta_{\e},\a)+o(1),
$$
where $o(1)$ is a family tending to zero; here, and in what
follows, unless otherwise stated, all limits are taken as $\e
\rightarrow 0$.

(2) The second one corresponds to the total error probabilities. Let
$\g_{\e}(\Theta_{\e},\psi)$ be the sum of the type I and the maximal
type II error probabilities, and let $\g_{\e}(\Theta_{\e})$ be the
minimax total error probability, i.e.,
$$
\g_{\e}(\Theta_{\e}) = \inf_{\psi}\g_{\e}(\Theta_{\e},\psi),
$$
where the infimum is taken over all possible tests. We call a
family of tests $\psi_{\e}$ {\it asymptotically minimax} if
$$
\g_{\e}(\Theta_{\e},\psi_{\e})=\g_{\e}(\Theta_{\e})+o(1).
$$
It is known that (see, e.g., \cite{IS.02}, Chapter 2) that
\begin{equation}\label{bg}
\b_{\e}(\Theta_{\e},\a)\in [0,1-\a], \quad
\g_{\e}(\Theta_{\e})=\inf_{\a\in (0,1)}(\a+\b_{\e}(\Theta_{\e},\a))
\in [0,1].
\end{equation}

We consider the problems of rate and sharp asymptotics for the error
probabilities in the minimax setup. The rate optimality problem
corresponds to the study of the conditions for which
$\g_{\e}(\Theta_{\e})\to 1$ and $\g_{\e}(\Theta_{\e})\to 0$ and,
under the conditions of the last relation, to the construction of
{\em asymptotically minimax consistent}\, families of tests $\psi_{\e}$, i.e.,
such that $\g_{\e}(\Theta_{\e},\psi_{\e})\to 0$.

We are interesting in a set $\Theta_{\e}$ of the form
$$
\Theta_\e=\Theta(r_{\e})=\{\eta\in \Theta\ : |\eta|\ge r_\e\},
$$
where $\Theta\subset l_2$ is a given set, $|\cdot |$ is some norm in
$l_2$ (not necessarily the standard $l_2$-norm) and $r_\e\to 0$ is a
given positive-valued family. For this case, we use the notation
$\g_\e(\Theta(r_\e))=\g_\e(r_\e),\
\b_\e(\Theta(r_{\e}),\a)=\b_\e(r_{\e},\a)$ and we are interesting in
the minimal decreasing rates for the sequence $r_{\e}$ such that
$\g_\e(r_\e)\to 0$. Namely, we say that the positive sequence
$r_{\e}^* \rightarrow 0$ is a {\it separation rate}, if
\begin{equation}\label{rates1}
 \g_{\e}(r_{\e})\to 1, \quad\text{and}\ \ \b_{\e}(r_{\e},\a)\to 1-\a \ \
\text{for any}\ \a\in (0,1),\quad\text{as}\quad r_{\e}/r_{\e}^*\to
0,
\end{equation}
and
\begin{equation}\label{rates2}
\g_{\e}(r_{\e})\to 0, \quad\text{and}\ \ \b_{\e}(r_{\e},\a)\to 0 \ \
\text{for any}\ \a\in (0,1),\quad \text{as}\quad r_{\e}/r_{\e}^*\to
\infty.
\end{equation}
In other words, it means that, for small $\e$, one can detect all
sequences $\eta\in \Theta(r_{\e})$ if the ratio $r_{\e}/r_{\e}^*$ is
large, whereas if this ratio is small then it is impossible to
distinguish between the null and the alternative hypothesis, with
small minimax total error probability. Hence, the rate optimality
problem corresponds to finding the separation rates $r_{\e}^*$ and
to constructing asymptotically minimax consistent families of tests.

On the other hand, the sharp optimality problem corresponds to the
study of the asymptotics of the quantities
$\b_{\e}(\Theta_{\e},\a),\ \g_{\e}(\Theta_{\e})$ (up to vanishing
terms) and to the construction of asymptotically minimax families of
tests $\psi_{\e,\a}$ and $\psi_{\e}$, respectively. Often, the sharp
asymptotics are of Gaussian type, i.e.,
\begin{equation}\label{G}
\b_{\e}(\Theta_{\e},\a)=\Phi(H^{(\a)}-u_{\e})+o(1),\quad
\g_{\e}(\Theta_{\e})=2\Phi(-u_{\e}/2)+o(1),
\end{equation}
where $\Phi$ is the standard Gaussian distribution function,
$H^{(\a)}$ is its $(1-\a)$-quantile, i.e.,
$\Phi(H^{(\a)})=1-\alpha$. The quantity $u_\e=u_\e(r_\e)$ is the value of the specific extreme
problem \nref{D1} on the sequence space $l^2$, and the extreme
sequence of this problem determines the structure of the asymptotically
minimax families of tests $\psi_{\e,\a}$ and $\psi_{\e}$.  Moreover,
we shall see that if $u_\e(r_\e)\to\infty$, then $\g_\e(r_\e)\to 0,\
\b_\e(r_\e,\a)\to 0$, and if $u_\e(r_\e)\to 0$, then $\g_\e(r_\e)\to
1,\ \b_\e(r_\e,\a)\to 1-\a$, for any $\a\in (0,1)$, i.e., the family
$u_{\e}(r_{\e})$ characterizes {\em distinguishability} in the
testing problem. The separation rates $r_{\e}^*$ are usually
determined by the relation $u_{\e}(r_{\e}^*)\asymp 1$ (see, e.g.,
\cite{I.93}, \cite{IS.02}). Hence, sharp and rate optimality
problems correspond to the study of the extreme problem \nref{D1}
and of the asymptotics of  the family $u_{\e}(r_\e)$.

\section{Minimax image detection from noisy tomographic data} \label{secnonadF}

\subsection{A general result: rate and sharp asymptotics} \label{subsec:genresF}
Recall the Gaussian sequence model \nref{ffgsm}. We are interested in the hypothesis testing
problem \nref{1.4} with the set under the alternative
$\Theta_\e=\Theta(r_\e)$ given by \nref{ell}.

Consider now the extreme problem
\begin{equation}\label{D1}
u_\varepsilon^2=u_\e^2(r_\e)=\frac{1}{2\e^4}\inf_{\eta\in\Theta(r_\e)}\sum_{\nu\in\Gamma}\eta_{\nu}^4.
\end{equation}
%
%
Suppose that $\Theta(r_\e)\not=\emptyset$ and $u_\e>0$, and let
there exist an extreme sequence $\{\tilde{\eta}_{\nu}\}_{\nu\in\Gamma}$ in the extreme problem
\nref{D1}. (Observe the uniqueness of a nonnegative extreme sequence
$\{\tilde{\eta}_{\nu}\}$, ${\nu\in\Gamma}$, because, by passing to the sequence $\{z_{\nu}\}_{\nu\in\Gamma}$  with elements $z_\nu=\tilde{\eta}_{\nu}^2$, ${\nu\in\Gamma}$, we obtain the minimization problem of a strictly convex function under linear constraints.) Denote
\begin{equation}\label{seq1}
w_{\nu}=\frac{\tilde{\eta}_{\nu}^2}{\sqrt{2\sum_{\nu\in\Gamma}\tilde{\eta}_{\nu}^4}},\,\, \nu\in\Gamma,
\quad w_0=\sup_{\nu\in\Gamma} w_{\nu},
\end{equation}
and consider the following families of test statistics and tests
\begin{equation}\label{test1}
t_\e=\sum_{\nu\in\Gamma}w_{\nu}\l((y_{\nu}/\e)^2-1\r),\quad \psi_{\e,H}=\1_{\{t_\e>H\}},
\end{equation}
where $\1_{\{A\}}$ denotes  the indicator function of a set $A$. (Note that the values
of $\tilde \eta_{\nu}$, $w_{\nu}$, ${\nu} \in \Gamma$, and $w_0$ depend on $\e$, i.e.,
$\tilde \eta_{\nu}= \tilde \eta_{{\nu},\e}$, $w_{\nu}=w_{{\nu},\e}$, ${\nu} \in \Gamma$, and $w_0=w_{0,\e}$.)

\medskip
The key tool for the study of the above mentioned hypothesis testing problem is the
following general theorem. Its proof follows along the lines
of the proof of Theorem 4.1 in \cite{ISS}; hence, it is omitted.

\begin{theorem}\label{testing}

Consider the Gaussian sequence model \nref{ffgsm} and the hypothesis
testing problem \nref{1.4} with the set under the alternative given
by \nref{ell}. Let $u_{\e}$ be determined by the extreme problem
\nref{D1}, let the coefficients $w_{\nu}$, $\nu\in\Gamma$, and $w_0$ be as in
\nref{seq1}, and consider the family tests $\psi_{\e,H}$ given by
\nref{test1}.  Then

(1) (a) If $u_\e\to 0$, then $\b_\e(r_\e,\a)\to 1-\a$ for any $\a\in
(0,1)$ and $\g_\e(r_\e)\to 1$, i.e., minimax testing is impossible.
If $u_\e=O(1)$, then $\liminf\b_\e(r_\e,\a)>0$ for any $\a\in (0,1)$
and $\liminf\g_\e(r_\e)>0$, i.e., minimax consistent testing is
impossible.

(b) If $u_\e\asymp 1$ and $w_0=o(1)$, then the family of tests
$\psi_{\e, H}$ of the form \nref{test1} with $H=H^{(\a)}$ and
$H=u_\e/2$ are asymptotically minimax, i.e.,
\begin{eqnarray*}
\a_{\e}(\psi_{\e,H^{(\a)}}) &\leq&  \a+o(1),\\
\b_{\e}(\Theta(r_{\e}), \psi_{\e,H^{(\a)}})&=&\b_{\e}(r_{\e},\a)+o(1),\\
\g_{\e}(\Theta(r_{\e}), \psi_{\e,u_\e/2})&=&\g_{\e}(r_{\e})+o(1),
\end{eqnarray*}
and the sharp asymptotics \nref{G} hold true, i.e.,
\begin{eqnarray*}
\b_{\e}(r_{\e},\a)&=&\Phi(H^{(\a)}-u_\e)+o(1), \\
\g_\e(r_{\e})&=&2\Phi(-u_\e/2)+o(1).
\end{eqnarray*}

(2) If $u_\e\to\infty$, then the family of tests $\psi_{\e,H}$ of
the form \nref{test1} with $H=T_\e$ are asymptotically minimax
consistent for any $c\in (0,1)$ and a family $T_\e\sim cu_\e$, i.e.,
$\g_\e(\Theta(r_{\e}),\psi_{\e,T_\e})\to 0$.
\end{theorem}

Theorem \ref{testing} shows that the asymptotics of the quality of
testing is determined by the asymptotics of values $u_\e$ of the the
extreme problem \nref{D1}. In order to make use of it, one needs to study
the extreme problem \nref{D1}. This problem is studied by using
Lagrange multipliers. Then, the extreme sequence in the above
mentioned extreme problem is of the form
\begin{equation}\label{seq}
\tilde{\eta}_{\nu}^2=z_0^2\sigma_{\nu}^2(1-Aa_{\nu}^2)_+, \quad \nu\in\Gamma,
\end{equation}
where $(t)_+=\max(t,0),\ t\in\R$, and  the quantities $z_0=z_{0,\e}$
and $A=A_{\e}$ are determined by the equations
\begin{equation}\label{eq0}
\begin{cases} &\sum_{\nu\in\Gamma} \sigma_{\nu}^2 \tilde{\eta}_{\nu}^2=r_\e^2,\cr
&\sum_{\nu\in\Gamma} a_{\nu}^2\sigma_{\nu}^2 \tilde{\eta}_{\nu}^2=1.
\end{cases}
\end{equation}
The equations \nref{eq0} are immediately rewritten in the form
\begin{equation}\label{eq}
\begin{cases} &r_\e^2=z_0^2J_1,\cr
&1=z_0^2A^{-1}J_2,
\end{cases}
\end{equation}
and, hence, the extreme problem \nref{D1} takes the form
\begin{equation}\label{eq1}
u_\varepsilon^2=\e^{-4}z_0^4J_0/2,
\end{equation}
where
\begin{eqnarray*}
J_1 &=& \sum_{\nu\in\Gamma}\sigma_{\nu}^4(1-A a_{\nu}^2)_+, \\
J_2 &=& A\sum_{\nu\in\Gamma}a_{\nu}^2\sigma_{\nu}^4(1-A a_{\nu}^2)_+,\\
J_0 &=& J_1-J_2=\sum_{\nu\in\Gamma}\sigma_{\nu}^4(1-A a_{\nu}^2)_+^2.
\end{eqnarray*}
It is also convenient to rewrite \nref{eq} and \nref{eq1} in the form
\begin{equation}
\label{eqJA}
r_{\e}^2 =A\, \frac{J_1}{J_2}, \quad u_{\e}^2 = \left(\frac{r_{\e}}{\e}\right)^4
\frac{J_0}{2J_1^2}.
\end{equation}


\begin{remark}
\label{rescalargFF}
{\rm Let $u_{\e}=u_{\e}(r_{\e})$ be the value of the extreme problem \nref{D1}
 with sequences $a=\{a_{\nu}\}_{\nu\in\Gamma}$ and
$\sigma=\{\sigma_{\nu}\}_{\nu \in\Gamma}$ associated with the set under the alternative
$\Theta_\e=\Theta(r_\e)$ given by \nref{ell}, and let $\tilde u_{\e}=\tilde u_{\e}(r_{\e})$
 be the corresponding value of the extreme problem similar to \nref{D1} with sequences
 $\tilde a=C a=\{C a_{\nu}\}_{\nu\in\Gamma}$ and
 $\tilde \sigma=D \sigma=\{D\sigma_{\nu}\}_{\nu\in\Gamma}$ in \nref{ell}, for some  constants $C, D>0$.
 Then, it is easily seen that the relation
$
\tilde u_{\e}(r_{\e}) = (CD)^{-2}u_{\e}(C r_{\e})
$
holds true.}
\end{remark}

\begin{remark}
\label{orderFF}
{\rm In order to obtain the corresponding rate and sharp asymptotics for the noisy tomographic data, we need to study the asymptotics of the quantities $J_i$, $i=0,1,2$, given above. We note, however, that the methods used in \cite{ISS} to study analogous asymptotics in a wide-range of linear statistical ill-posed inverse problems cannot be adopted to the problem at hand. The reason is that there does not exist a common ordering for the sequences  $a=\{a_{\nu}\}_{\nu\in\Gamma}$ and
$\sigma=\{\sigma_{\nu}\}_{\nu \in\Gamma}$ associated with the set under the alternative
$\Theta_\e=\Theta(r_\e)$ given by \nref{ell}. The arguments and techniques used to prove Theorems \ref{T1} and \ref{T2} below are specifically developed to tackle this problem.
}
\end{remark}


\subsection{Rate and sharp asymptotics for the noisy tomographic data}
According to \nref{radsing} and \nref{functthlF}, consider the
double-index sequences
\begin{eqnarray}\label{tom2}
 a_{\nu} &=& L^{-1}(j+1)^{p}(l+1)^{p}, \,\, \,\, \nu \in \Gamma,\\
 \label{tom3}
\sigma_{\nu} &=& \pi^{-1}(j+l+1)^{1/2}, \quad  \,\,\, \nu \in
\Gamma,
 \end{eqnarray}
for some $p>0$ with $\Gamma$ given by \nref{GammaFF}.

\begin{theorem}\label{T1}
Consider the Gaussian sequence model \nref{ffgsm} and the hypothesis
testing problem \nref{1.4} with the set under the alternative given
by \nref{ell}. Let $\{a_{\nu}\}_{\nu \in \Gamma}$ and $\{\sigma_{\nu}\}_{\nu \in \Gamma}$ be defined as in \nref{tom2}  and \nref{tom3}, respectively. Then

(a) The sharp asymptotics \nref{G} hold with  the value $u_\e$ of
the extreme problem \nref{D1} determined by
\begin{equation}
\label{UMFF2} u_\e^2\sim \pi^4 L^{-3/p}
r_\e^{4+3/p}\e^{-4}\,\frac{2p+3}{2B}\left(\frac{3}{4p+3}\right)^{1+3/(2p)},
\end{equation}
where $B=\sum_{m=1}^{\infty}m^{-3},\ 1.2021<B<1.2022$.

(b) The asymptotically minimax family of tests $\psi_{\e,H}$ are
determined by the family of test statistics $t_{\e}$ given by
\nref{test1} with coefficients $w_{\nu}$, $\nu\in \Gamma$, and $w_0$
as in \nref{seq1}, and with extreme sequence
$\{\tilde{\eta}_{\nu}\}_{\nu\in\Gamma}$ satisfying \nref{seq} with
$ A \sim \,\frac{3}{4p+3}r_\e^2$.

(c) The separation rates are of the form
\begin{equation}\label{seprateSobSovF}
 r_\varepsilon^* = \varepsilon^{4p/(4p+3)}.
\end{equation}
\end{theorem}

\begin{remark}
\label{FFF1FF}
{\rm It is easily seen that the asymptotic results in Theorem \ref{T1} hold true uniformly over
$p \in \Sigma$ for any compact set $\Sigma\subset (0,\infty)$.
}
\end{remark}

\begin{remark}
\label{FFF1FFanis1}
{\rm Rate and sharp asymptotics in the corresponding minimax estimation problem under the $L^2$-risk have been obtained in \cite{Cav-Tsyb}. In particular, the asymptotical (as $\e \rightarrow 0$) minimax rates of estimation are given by
$$
R^2_{\e}:=\inf_{\tilde{f}} \sup_{f \in {\cal F}(p,L)} E||\tilde{f}-f||^2 \asymp\e^{4p/(2p+2)},
$$
where the infimum is taken over all possible estimators $\tilde{f}$ of $f$ based on observations from the
Gaussian white noise model \nref{gwmRDTFan}. (Here, we adopt standard notation and write $g_1(\e) \asymp g_2(\e)$
to denote $0 < \liminf(g_1(\e)/g_2(\e)) \leq \limsup(g_1(\e)/g_2(\e))
< \infty$ as $\e \rightarrow 0.$) By comparing $r_\varepsilon^*$ with $R_{\e}$, it is observed that the asymptotical minimax rates of testing are faster than the corresponding asymptotical minimax rates of estimation; this phenomenon is common in nonparametric statistical inference (see, e.g., \cite{IS.02}, Sections 2.10 and 3.5.1, \cite{ISS}).
}
\end{remark}

\subsection{Adaptivity and rate optimality for the noisy tomographic data}
The family of tests considered in Section \ref{T1} depends on the parameter $p$ that is usually
unknown in practice. Therefore, it is of paramount importance to construct families of tests
that do not depend on the unknown parameter $p$ and, at the same time, provide the best possible
asymptotical minimax efficiency. Such families are called {\em adaptive} (to the parameter $p$) and
the formal setting is as follows.

Let $\Sigma$ be  a compact set in
$(0,\infty)$ and a family $r_\e(p),\ p \in\Sigma$, be given, where
$\e>0$ is small. Let the set $\Theta_{\e}(p, r_{\e}(p))$ be
determined by the constraints \nref{ell} with $a_{\nu}=a_{\nu}(p), \
\nu \in\Gamma$, and $r_\e=r_{\e}(p)$, and set
$$
\Theta_{\e}(\Sigma)=\bigcup_{p\in\Sigma}\Theta_\e(p,r_{\e}(p)).
$$
We are interesting in the following hypothesis testing problem
$$
H_0:\ \eta=0,\quad \text{versus}\quad H_1:\ \eta\in
\Theta_{\e}(\Sigma).
$$
We aim at finding conditions for either
$\g_\e(\Theta_{\e}(\Sigma))\to 1$ or $\g_\e(\Theta_{\e}(\Sigma))\to
0$, and  to constructing asymptotically minimax adaptive consistent families of tests
$\psi^{ad}_\e$ such that
$\g_\e(\Theta_{\e}(\Sigma),\psi^{ad}_\e)\to 0$ as
$\g_\e(\Theta_{\e}(\Sigma))\to 0$.

Let $u_\e(p)=u_\e(p, r_{\e}(p))$ be the value of the extreme problem
\nref{D1} for
the set $\Theta_\e=\Theta_{\e}(p,r_{\e}(p))$. Set
$$
u_\e(\Sigma)=\inf_{p\in\Sigma}u_\e(p).
$$
We are interesting in how large $u_\e(\Sigma)$ should be in order to
provide the relation $\g_\e(\Theta_{\e}(\Sigma))\to 0$.
We say that the family $u_\e^{ad}=u_\e^{ad}(\Sigma)\to\infty$
characterizes {\it adaptive distinguishability} if there exist
constants $0<d=d(\Sigma)\le D(\Sigma)=D<\infty$ such that
\begin{eqnarray*}
\g_\e(\Theta_{\e}(\Sigma))\to 1 & \text{as}\quad
\lim\sup_{p\in\Sigma}u_\e(p)/u_\e^{ad}<d,\\
\g_\e(\Theta_{\e}(\Sigma))\to 0 & \text{as}\quad
\lim\inf_{p\in\Sigma}u_\e(p)/u_\e^{ad}>D.
\end{eqnarray*}

Observe that it follows from the asymptotics \nref{UMFF2} that, by
making $r_\e(p)$ larger or smaller, one can increase or decrease
$u_\e(p,r_\e(p))$ in order to get $u_\e(p,r_\e(p))\sim u_\e$, for
all $p\in \Sigma$ and any family $u_\e>0$.

We call  a family $r^{ad}_\e(p),\ p\in\Sigma,$ such that
$u_\e^{ad}\asymp u_\e(p,r^{ad}_\e(p))$, the family of {\it adaptive
separation rates}.

Note that the relation $\g_\e(\Theta_{\e}(\Sigma))\to 0$ is possible
if $u_\e(\Sigma)\to\infty$. However this implication does
not hold for the tomography problem under consideration, as we show below.
(A similar situation appears in some ill-posed inverse problems, see \cite{ISS}).) Hence,
hereafter,  adaptive distinguishability
conditions and adaptive separation rates are sought for the tomography problem. In contrast
 to Theorem \ref{T1}, there is price to pay for the adaptation. We show below that
\begin{equation}
\label{adcFF}
u_\e^{ad}=\sqrt{\log\log\e^{-1}},
\end{equation}
yielding a loss in the separation rates in terms of an extra factor
$\sqrt[4]{\log\log\e^{-1}}$ in $\e$. Furthermore, the derived
families of tests are of simple structure.
 (A similar loss in the separation rates was first
observed in \cite{Spok} and more recently in \cite{ISS}.)

Specifically, let $p$ be unknown, $p\in \Sigma$, where
$\Sigma=[p_{min},p_{max}],\ 0<p_{min}<p_{max}<\infty$, be a compact
interval in $(0,\infty)$. Let us also consider the collections
$$
p_k\in\Sigma,\quad c_k\sim 2/r_\e(p_k),\quad k=0,1,\ldots,K, \quad
K=K_\e\asymp \log(\e^{-1})\log\log(\e^{-1}),
$$
where $p_0=p_{max}>p_1>\ldots>p_K=p_{min}$ (the collection $p_k,\
k=1,2,\ldots, K-1$, will be specified in the proof)
 and collection of statistics $t_{\e,c_k}$ of the form
\begin{equation}\label{adtom5}
t_{\e,c_k}=\sum_{\nu \in C_{\nu,k}}w_{\nu,k}((y_\nu^2/\e^2)-1),\,\,
w_{\nu,k}=\frac{\sigma_\nu^2}{\left(2\sum_{\nu \in C_{\nu,k}}\sigma_\nu^4\right)^{1/2}},\,\,
\sum_{\nu \in C_{\nu,k}} w_{\nu,k}^2=\frac12,
\end{equation}
for $\nu \in \Gamma$ given by  \nref{GammaFF} and
$C_{\nu,k}=\{\nu:\,a_{\nu,p_k}\le c_k\}$. Consider the following
families of thresholds and tests
\begin{equation}\label{adtomF}
H_\e=2\sqrt{\log(K_\e)},\quad \CY_\e=\{y:\ t_{\e, c_k}\le H_\e,\quad
\forall\ 0\le k\le K_\e\},\quad \psi_\e=\1_{{\overline \CY}_\e},
\end{equation}
where $\bar A$ denotes the complement of a set $A$.

Denote also
\begin{equation}
\label{phiFF}
\phi(p)=\frac{4}{4p+3},\quad \phi(\Sigma)=\{\phi(p):\,
p\in\Sigma\}\subset (0,\infty).
\end{equation}

\begin{theorem}\label{T2}
Consider the Gaussian sequence model \nref{ffgsm} and the hypothesis
testing problem \nref{1.4} with the set under the alternative given
by \nref{ell}. Let $\{a_{\nu}\}_{\nu \in \Gamma}$ and $\{\sigma_{\nu}\}_{\nu \in \Gamma}$ be defined as in \nref{tom2}  and \nref{tom3}, respectively. Then

(a) (lower bounds) Let the set $\phi(\Sigma)$ given by \nref{phiFF} contains an interval $[a,b],\ 0<a<b<4/3$. There exists constant $d=d(\Sigma)>0$ such that if $\limsup
_{p\in\Sigma}u_\e(p)/\sqrt{\log\log(\e^{-1})}\le d$, then
$\g_\e(\Theta_\e(\Sigma))\to 1$.

(b) (upper bounds) For the family of tests $\psi_{\e}$ given by
\nref{adtomF}, $\a(\psi_\e)=o(1)$ and there exists constant
$D=D(\Sigma)>0$ such that if
$\liminf_{p\in\Sigma}u_\e(p)/\sqrt{\log\log(\e^{-1})}>D$, then
$\b_\e(\Theta_\e(\Sigma),\psi_\e)=o(1)$.

(c) (adaptive separation rates) The adaptive distinguishability family $u_\e^{ad}$ is given by \nref{adcFF} and the adaptive separation rates $r_{\e}^{ad}(p)$,
$p \in \Sigma$, are
given by
$$
r^{ad}_\e(p)=\left(\e\sqrt
[4]{\log\log(\e^{-1})}\right)^{4p/(4p+3)}.
$$
\end{theorem}


\begin{remark}
\label{FFF1FFanis2}
{\rm Rate and sharp adaptation in the corresponding minimax estimation problem under the $L^2$-risk have been obtained in \cite{Cav-Tsyb}. In particular, \cite{Cav-Tsyb} showed that
$$
R^2_{\e}=\rho_{\e}(p,L) (1+o(1)),
$$
where
$$
\rho_{\e}(p,L)= \frac{1}{2} \bigg(\frac{\pi^4 p}{3(p+2)}\bigg)^{2p/(2p+2)} ((2p+2)L)^{2/(2p+2)}\, \e^{4p/(2p+2)},
$$
and constructed an adaptive penalized blockwise Stein-type estimator $\hat{f}$ (see \cite{Cav-Tsyb}, equations (3.12) and (5.8)) such that
$$
\lim_{\e \rightarrow 0}\,\, \sup_{p \in [p_1,p_2], \, L \in [L_1,L_2]} \,\, \sup_{f \in {\cal F}(p,L)}
\frac{E||\hat{f}-f||^2}{\rho_{\e}(p,L)}=1,
$$
for any $0 <p_1<p_2<\infty$ and $0<L_1<L_2<\infty$. Note that, in contrast to adaptive minimax separation rates, there is no price to pay for adaptation in the corresponding minimax estimation problem under the $L^2$-risk (i.e., a global measure). The unavoidable logarithmic factor for adaptivity that appears in the minimax separation rates $r_{\e}^{ad}(p)$ stated in Theorem \ref{T2} also appears in various other adaptivity problems, such as minimax signal detection (see \cite{ISS}) and minimax estimation under the $L^2$-risk (see \cite{T2000}) in some ill-posed inverse problems. It also resembles the minimal price one needs to pay for adaptation in minimax
estimation under the $l^2$-risk (i.e., a local or pointwise measure) that has been observed in  \cite{L1990}, \cite{br-lo}, \cite{LS1997} and \cite{T1998} in the case of Lipschitz and Sobolev balls and, more recently, in \cite{Cai2003} and \cite{BS2009} in the case of Besov balls. }
\end{remark}

\section{Appendix: Proofs}

For simplicity in the calculations, we omit the factors $L^{-1}$ and
$\pi^{-1}$ in \nref{tom2} and \nref{tom3}, respectively. In other
words, from now onwards, we work with $a_{\nu}=(j+1)^{p}(l+1)^{p}$
and $\sigma_{\nu}=(j+l+1)^{1/2}$, $p>0$, $\nu \in \Gamma$. The final
results can be obtained on rescaling by using Remark
\ref{rescalargFF}.

\subsection{\em Proof of Theorem \ref{T1}}\label{S1}

It follows from Theorem \ref{testing} that the efficiency in the detection problem under consideration is
determined by the asymptotics of the quantity
$$
u_\varepsilon^2=\left(\frac{r_\varepsilon}{\varepsilon}\right)^4\frac{J_0}{2J_1^2},
$$
where
\begin{eqnarray*}
J_1=J_1(A)&=&\sum_{\nu \in \Gamma}\sigma_{\nu}^4(1-Aa_{\nu}^2)_+,\label{tom4}\\
J_2=J_2(A)&=&A\sum_{\nu \in \Gamma} a_\nu^2\sigma_\nu^4(1-Aa_\nu^2)_+,\label{tom5}\\
J_0=J_0(A)&=&\sum_{\nu \in \Gamma} \sigma_\nu^4(1-Aa_\nu^2)^2_+=J_1-J_2,\label{tom6}
\end{eqnarray*}
for $\nu \in \Gamma$ given by  \nref{GammaFF}). Moreover, the quantity $A=A_\e\to 0$ is
determined by the relation
\begin{equation}\label{tom8}
r_\varepsilon^2=A\frac{J_1}{J_2}.
 \end{equation}
 In order to study the asymptotics of $u_\e$,
we are interesting in the asymptotics of the functions $J_i(A)$, $i=0,1,2$, as $A\to 0$.
We first, however, start with the asymptotics of the following function
\begin{equation}
I(A)=\sum_{\{\nu:\ Aa_\nu^2\leq 1\}}\sigma_\nu^4, \quad \nu \in \Gamma.\label{tom7}
 \end{equation}

\begin{proposition}\label{P1}
Let $I(A)$ be defined as in \nref{tom7}. Then, as $A\to 0$,
\begin{equation}\label{tom8a}
 I(A)\sim\frac{2B}3A^{-3/(2p)},
\end{equation}
where  $B=\sum_{m=1}^{\infty}m^{-3},\ 1.2020<B<1.2021$.
\end{proposition}
{\bf Proof}. Set $j+1=m$, $l+1=n$, $H=[m^{-1}A^{-1/(2p)}]$ and
$H_1=[A^{-1/(4p)}]$, where $[t]$ is the integer part of  $t$.
Consider the set
$$
C_{m,n,p,A}=\{(m,n):\,m\ge 1, n\ge 1,(mn)^{2p}\leq
A^{-1}\}.
$$
Then, we have
\begin{eqnarray*}
I(A)&=&\sum_{(m,n)\in C_{m,n,p,A}}(m+n-1)^2\\
&=&2\sum_{m=1}^{H_1}\sum_{n=1}^{H}(m+n-1)^2-\sum_{m=1}^{H_1}\sum_{n=1}^{H_1}(m+n-1)^2\\
&=&\frac{1}{3}\sum_{m=1}^{H_1}((m-1+H)(m+H)(2m-1+2H)-(m-1)m(2m-1))+O(A^{-1/p})\\
&=&2\sum_{m=1}^{H_1}(Hm^2+H^2m-H m -H^2/2+H^3/3+H/6)+O(A^{-1/p}).
\end{eqnarray*}

We now indicate the asymptotics of the items in the last sum. Observe that
$$
H=\frac{1}{m}A^{-1/(2p)}+\alpha(A,m),\quad H_1=A^{-1/(4p)}+\beta(A),
$$
where $ \alpha(A,m)\in [0,1)$ and $\beta(A)\in [0,1). $
Thus, we have
\begin{eqnarray*}
\sum_{m=1}^{H_1}
Hm^2&=&\sum_{m=1}^{H_1}(A^{-1/(2p)}m+\alpha(A,m)m^2)=A^{-1/p}/2+O(A^{-3/(4p)}),\\
\sum_{m=1}^{H_1}
H^2m&=&\sum_{m=1}^{H_1}(A^{-1/(2p)}/m+\alpha(A,m))^2\,m\sim
A^{-1/p}\log(A^{-1/(4p)})=
\frac{A^{-1/p}\log(A^{-1})}{4p},\\
\sum_{m=1}^{H_1} Hm&=&\sum_{m=1}^{H_1}
(A^{-1/(2p)}+\alpha(A,m)m)\sim A^{-3/(4p)},\\
\sum_{m=1}^{H_1} H^2/2&=&\frac12\sum_{m=1}^{H_1}
(A^{-1/(2p)}/m+\alpha(A,m))^2\asymp A^{-1/p},\\
\sum_{m=1}^{H_1} H^3/3&=&\frac{A^{-3/(2p)}}{3}\sum_{m=1}^{H_1}
m^{-3}+O(A^{-1/p}),\\
\sum_{m=1}^{H_1} H/6&=&\sum_{m=1}^{H_1}
(A^{-1/(2p)}/m+\alpha(A,m))/6=o(A^{-1/p}).
\end{eqnarray*}
Therefore,
\begin{equation*}
I(A)=\frac{2A^{-3/(2p)}}{3}\sum_{m=1}^{H_1}m^{-3}+\frac{A^{-1/p}\log(A^{-1})}{2p}
+O(A^{-1/p}) \sim\frac{2B}3A^{-3/(2p)}.
\end{equation*}
The proposition now follows. \quad\endproof

\medskip
Let us now return to the asymptotics of $J_i(A)$, $i=0,1,2$, as $A \to 0$. Introduce
the following function
\begin{equation}\label{tom8ba}
F(t)=I(t^{-1})=\sum_{(m,n)\in C_{m,n,p,t^{-1}}}(m+n-1)^2,\quad t\ge 0,
\end{equation}
and observe that  $F(t)$ is nondecreasing in $t\ge 0$ such that $F(0)=0$. It
follows from Proposition \ref{P1} that
\begin{equation}\label{tom8b}
F(t)\sim \frac{2B}3t^{3/(2p)},\quad t\to\infty.
\end{equation}
For $T=A^{-1}$, the functions $J_i(A)$, $i=0,1,2$, could be rewritten in the
form
\begin{eqnarray*}
J_1(A)&=&\int_0^T(1-t/T)dF(t),\\
J_2(A)&=&\int_0^T(t/T-(t/T)^2)dF(t),\\
J_0(A)&=&\int_0^T(1-t/T)^2dF(t)=J_1-J_2.
\end{eqnarray*}
Integrating now by parts, we get
\begin{equation*}
J_1(A)=T^{-1}\int_0^{T}F(t)dt.
\end{equation*}
In order to study the asymptotics of the above integral, we divide it into
the following two parts
$$
\int_0^{T}F(t)dt=S_1+S_2,\quad S_1=\int_0^{T^{1/2}}F(t)dt,\quad
S_2=\int_{T^{1/2}}^{T}F(t)dt.
$$
Hence, it suffices to check that  $S_1=o(S_2)$ and to use the
asymptotics \nref{tom8b} of the function $F$ under the integral in
$S_2$. It is easily seen that
\begin{equation*} 0\le S_1 \le F(T^{1/2})T^{1/2}\asymp
T^{3/(4p)+1/2}
\end{equation*}
and that
\begin{equation*}
S_2 \sim \int_{T^{1/2}}^{T}\left(\frac{2B}3t^{3/(2p)})\right)dt\sim
\frac{4Bp}{3(2p+3)}T^{1+3/(2p)}.
\end{equation*}
Therefore, we get
\begin{equation}\label{tom9}
J_1(A)\sim\frac{4Bp}{3(2p+3)}A^{-3/(2p)}.
\end{equation}
Similarly, for $T=A^{-1}$, we get
\begin{eqnarray}
J_2(A)&=&T^{-1}\int_0^{T}t dF(t) -T^{-2}\int_0^{T}t^2 dF(t)
=2T^{-2}\int_0^{T}tF(t)dt\nonumber\\
&-&T^{-1}\int_0^{T}F(t)dt\nonumber
\sim T^{3/(2p)}\left(\frac{8Bp}{3(4p+3)}-\frac{4Bp}{3(2p+3)}\right)\\
&=&\frac{4Bp}{(4p+3)(2p+3)}A^{-3/(2p)}\label{tom9a}
\end{eqnarray}
and that
\begin{equation}
J_0 = J_1-J_2\sim\frac{16Bp^2}{3(2p+3)(4p+3)}A^{-3/(2p)}\label{tom9b}.
\end{equation}
Thus, it follows from \nref{tom8}, \nref{tom9}, \nref{tom9a} and \nref{tom9b}
that
\begin{equation}\label{tom10}
A\sim \frac{3}{4p+3}r_\e^2,\quad u_\e^2\sim
r_\e^{4+3/p}\e^{-4}\frac{2p+3}{2B}\left(\frac{3}{4p+3}\right)^{1+3/(2p)}.
\end{equation}
Therefore, the separation rates  $r_\e^*$ are of the form
\begin{equation}\label{tom11}
r_\e^*=\e^{4p/(4p+3)}.
\end{equation}

In order to get the sharp asymptotics, in view of Theorem \ref{testing},
it is enough to check the condition
$$
w_0=\frac{\max_{\{\nu:\, a_\nu^2A<1\}}\sigma_\nu^2(1-Aa_\nu^2)}
{\sqrt{2\sum_{\{\nu:\, a_\nu^2A<1\}}\sigma_\nu^4(1-Aa_\nu^2)^2}}=o(1), \quad \nu \in \Gamma.
$$
This condition follows directly from the relation
\begin{equation}\label{tom12}
\max_{\{\nu:\, a_\nu^2A<1\}}\sigma_\nu^4=o(J_0(A)).
\end{equation}
Indeed, for $m=j+1\in \N,\,n=l+1\in\N$, the last condition follows from
(see \nref{tom9b})
$$
\sigma^4_\nu=(m+n-1)^2<4A^{-1/p}\ll A^{-3/2p}\asymp J_0(A)
$$
as  $A\to 0$, for $\nu=(m-1,n-1)$ such that $\ a_\nu^2=(mn)^{2p}\leq A^{-1}$.

The theorem now follows.

\subsection{\em Proof of Theorem \ref{T2}}\label{S11}
Let $p$ be unknown and consider $p\in \Sigma$, where
$\Sigma=[p_{min},p_{max}],\ 0<p_{min}<p_{max}<\infty$, is a compact
interval in $(0,\infty)$.  Recall the expression $\phi(p)$ from
\nref{phiFF} and that $Z_+=\N \cup \{0\}$.

We first obtain the lower bounds. Take a collection $p_k$, $k \in Z_+$, such that
$$
\phi(p_k)=a+k\delta_\e,\quad k=0,1,\ldots,K=K_{\e}
$$
with
$$
\phi(p_K)=b>a=\phi(p_0)>0,\quad
\delta=\delta_\e=\frac{(b-a)}{K}\sim\frac{\log(2)}{\log(\e^{-1})}.
$$
Therefore,
$$
p_k\in[b^{-1}-3/4,a^{-1}-3/4]=\Sigma,\quad b<4/3,
$$
and
$$
\delta_\e=\phi(p_k)-\phi(p_{k-1}), \quad
p_{k-1}-p_k=(4p_k+3)(4p_{k-1}+3)\delta_\e/16\asymp\frac1{\log(\e^{-1})}=o(1).
$$
Assume without loss of generality that, uniformly in $p\in \Sigma$,
\begin{equation}\label{M*2}
u_\e(p)\sim\sqrt{d\log\log(\e^{-1})},
\end{equation}
where the constant $d>0$ will be specified below. This corresponds
to taking, uniformly in $p\in \Sigma$,
\begin{equation}\label{M*1}
r_\e(p)\sim \l(\e\l(d(p)\log\log(\e^{-1})\r)^{1/4}\r)^{\phi(p)p},
\end{equation}
where
$$d(p)=da(p),\quad
a(p)=\frac{2B}{2p+3}\left(\frac34\phi(p)\right)^{-(2p+3)/(2p)}.
$$
Take
\begin{equation*}
T_k\sim (2r_\e(p_k))^{-1/p_k},\quad k=0,1,\ldots,K_\e.
\end{equation*}
By construction, we have
\begin{eqnarray*}
T_k-T_{k-1}&\sim&
T_{k-1}\l(2^{p_{k-1}^{-1}-p_k^{-1}}(\e\sqrt [4]{d\log\log(\e^{-1})})^{-\delta_\e}-1\r)\\
&=& T_{k-1}(\exp(\log(2)(1+o(1)))-1) \sim T_{k-1}.
\end{eqnarray*}
Observe that the function $F(t)=F_p(t)$ defined by \nref{tom8ba}
depends on $p$. Also, by \nref{tom8b}, we can write
\begin{eqnarray}\label{Fnew}
F(t)=F_p(t)&=&\sum_{\{((i_1+1)(i_2+1))^{2p}\le t,\ i_1,i_2\in
Z_+\}}(i_1+i_2+1)^2 \nonumber\\
&=&\sum_{\{\nu:\ a_\nu^2\le
t\}}\sigma_\nu^4=F_1(t^{1/p})\sim\frac{2B}{3}t^{3/(2p)}.
\end{eqnarray}

Set
$$
\Delta_k=\{\nu_0=(j,l)\in\Z_+\times\Z_+: T_{k-1}<(j+1)(l+1)\le
T_k\}.
$$
Take  a collection $z_k>0,\ v_{k,\nu_0},\ \nu_0\in\Delta_k$,
$k=1,2,\ldots,K$, such that
$$
v_{k,\nu_0}=z_k
\begin{cases}
&\xi_{\nu_0}\sigma_{\nu_0},\quad \text{if} \quad \nu_0\in
\Delta_k,\cr &0, \qquad \quad\text{otherwise},
\end{cases}
$$
with $\xi_{\nu_0}=\pm 1$, $\nu_0 \in \Delta_k$. We then have
\begin{equation}\label{constr_a}
\sum_{\nu_0 \in
\Delta_k}v_{k,\nu_0}^2\sigma_{\nu_0}^2=z_k^2\sum_{\nu_0\in\Delta_k}\sigma_{\nu_0}^4=
z_k^2\l(F_1(T_k^{ 2})-F_1(T_{k-1}^{2})\r)
\sim\frac{7B}{12}z_k^2T_{k}^3=2r^2_\e(p_k)
\end{equation}
and
\begin{equation}\label{constr_b}
\sum_{\nu_0 \in
\Delta_k}v_{k,\nu_0}^2\sigma_{\nu_0}^2a_{\nu_0,p_k}^2=z_k^2\sum_{\nu_0\in\Delta_k
}\sigma_{\nu_0}^4a_{\nu_0,p_k}^2\leq
2T_k^{2p_k}r^2_\e(p_k)(1+o(1))=\frac12(1+o(1))<1.
\end{equation}
Furthermore, we have

$$
u_\e^2(p_k)\sim d\,\log\log(\e^{-1}),\ d>0,\quad p_k\in\Sigma.
$$
Consider now the priors
$$
\pi_k=\prod_{\nu_0\in\Delta_k}(\delta_{z_k \sigma_{\nu_0}
e_{\nu_0}}+\delta_{-z_k \sigma_{\nu_0} e_{\nu_0}})/2, \quad
k=1,2,\ldots,K, \quad \pi=\frac{1}{K}\sum_{k=1}^K\pi_k,
$$
where $\{e_{\nu_0}\}_{\nu_0\in\Z_+\times\Z_+}$ is the standard basis in the
space $l^2$ that corresponds to sequences indexed by $\nu_0\in \Z_+\times\Z_+$, and $\delta_\eta$ is the Dirac mass at the
point $\eta\in l^2$. The relations \nref{constr_a} and \nref{constr_b}
imply
$$
\pi_k(\Theta_\e(p_k,r_\e(p_k)))=1,\quad \pi(\Theta_\e(\Sigma))=1.
$$
Let $P_{\pi_k}=E_{\pi_k}P_{\e,\eta},\ P_\pi=E_\pi P_{\e,\eta}$ be
the mixtures over the priors. It suffices to check that
\begin{equation}\label{l2norm}
E_{\e,0}\l(\l(dP_{\pi}/dP_{\e,0}-1\r)^2\r)=o(1).
\end{equation}
Using evaluations similar to Section 5.6 in \cite{ISS},  we have
\begin{eqnarray*}
E_{\e,0}\l(\l(dP_{\pi}/dP_{\e,0}-1\r)^2\r)&=&\frac{1}{K^2}
\sum_{k=0}^K E_{\e,0}\l(\l(dP_{\pi_k}/dP_{\e,0}-1\r)^2\r)\\
&=&\frac{1}{K^2}
\sum_{k=0}^K\l(E_{\e,0}\l(dP_{\pi_k}/dP_{\e,0}\r)^2-1\r) \\
&\le&\frac{1}{K^2} \sum_{k=0}^K\l(\exp\
\l(2\sum_{\nu_0\in\Delta_k}\sinh^2(z_k^2\sigma_{\nu_0}^2/2\varepsilon^2)\r)-1
\r).
\end{eqnarray*}
Note that $\sigma_{\nu_0}^2=j+l+1<T_k$ if ${\nu_0}\in\Delta_k$. Therefore,
uniformly over ${\nu_0}\in\Delta_k$, we have
$$
\frac{z_k^2\sigma_{\nu_0}^2}{\e^2}<\frac{z_k^2T_k}{\e^2}\asymp
r_\e(p_k)^{2+2/p_k}\e^{-2}
\asymp\e^{2/(4p_k+3)}(\log\log(\e^{-1}))^{(2p_k+2)/(4p_k+3)}=o(1),
$$
and (since $ \sinh^2(z_k^2\sigma_{\nu_0}^2/2\varepsilon^2)\sim
z_k^4\sigma_{\nu_0}^4/4\varepsilon^4$)
\begin{eqnarray*}
2\sum_{{\nu_0}\in\Delta_k}\sinh^2(z_k^2\sigma_{\nu_0}^2/2\varepsilon^2)&\sim&
z_k^2\e^{-4}r^2_\e(p_k)\sim\frac{24}{7B}2^{3/p_k}
r^{4+3/p_k}_\e(p_k)\e^{-4}\\
&\sim&\frac{24}{7B}2^{3/p_k}da(p_k)\log\log(\e^{-1}).
\end{eqnarray*}
One can take $d>0$ such that, for any $p_k\in\Sigma$,
$$
d\frac{24}{7B}2^{3/p_k}a(p_k)\le
d\frac{16}7\max_{k=0,1,...,K_\e}\left\{\frac{4p_k+3}{2p_k+3}
\left(\frac{4(4p_k+3)}3\right)^{3/(2p_k)} \right\}=d_1<1.
$$
Then we have
\begin{eqnarray*}
E_{\e,0}\l(\l(dP_{\pi}/dP_{\e,0}-1\r)^2\r) &\le&\frac{1}{K^2}
\sum_{k=0}^K\l(\exp\l(2\sum_{{\nu_0}\in\Delta_k}\sinh^2(z_k^2\sigma_{\nu_0}^2/2\varepsilon^2)\r)
-1\r)\\
&<&\frac{K\log^
{d_1}(\e^{-1})}{K^2}\asymp\frac{\log^{d_1}(\e^{-1})}{\log(\e^{-1})}=o(1).
\end{eqnarray*}

We now obtain the upper bounds.

Similarly to the proof of the lower bounds, assume without loss of
generality that $u_\e(p)\sim \sqrt{D\log\log(\e^{-1})}$, uniformly
in $p\in \Sigma$, where the constant $D>0$ will be specified below.
This corresponds to \nref{M*1} with $d$ replaced by $D$, uniformly
in $p\in \Sigma$.

In order to evaluate the type I error probability, we consider a
different grid with different $K=K_\e$, i.e.,
$$
\phi(p_k)=a+k\delta_\e,\quad k=0,1,\ldots,K=K_{\e},\quad \phi(p_K)=b>a>0,
$$
where
$$
\delta=\delta_\e=\frac{(b-a)}{K_\e}\sim\frac{\log(2)}{\log(\e^{-1})\log\log(\e^{-1})}.
$$

Let us evaluate the exponential moments
$$
E_{\e,0}(\exp(ht_{\e, c_k})),\ h>0.
$$
Recall that $(y_\nu/\e) \stackrel{iid}{\sim} \CN(0,1)$, $\nu \in
\Gamma$, under $P_0$. Recall the set $C_{\nu,k}=\{\nu:\
a_{\nu,p_k}\le c_k\}$, $\nu \in \Gamma$, $k=0,1,\ldots,K$. Let the
family $h=h_\e$ be taken in such way that
$$
h\max_{\nu\in C_{\nu,k}}w_{\nu,k}=o(1).
$$
Then, we have
\begin{eqnarray}\label{adtom1}
E_{\e,0}(\exp(ht_{\e, c_k}))&=&\prod_{\nu \in
C_{\nu,k}}\left(\exp(-hw_{\nu,k})E_{\e,0}\exp(hw_{\nu,k}\xi_\nu^2)\right)
\nonumber\\
&=& \exp\left(\sum_{\nu \in C_{\nu,k}}(-hw_{\nu,k}-\log(1-2hw_{\nu,k})/2)\right)\nonumber\\
&=&\exp\left(\sum_{\nu \in
C_{\nu,k}}h^2w_{\nu,k}^2(1+O(hw_{\nu,k}))\right)\nonumber\\
&=&\exp\l(h^2/2)(1+o(1)\r).\qquad
\end{eqnarray}
Let $h=H_\e$. Then, for $k = 0, 1,\ldots, K_\e$, we have
\begin{eqnarray*}
h\max_{\nu \in C_{\nu,k}}w_{\nu,k} &=&H_\e\frac{\max_{\nu \in
C_{\nu,k}} \sigma_\nu^2}{\left(2\sum_{\nu \in
C_{\nu,k}}\sigma_\nu^4\right)^{1/2}}<\frac{H_{\e}
c_k^{1/p_k}}{\left(4Bc_k^{3/p_k}(1+o(1))/3\right)^{1/2}}\\
&\asymp& H_{\e}r^{1/(2p_k)}_\e(p_k)\asymp
\e^{2/(3+4p_k)}(\log\log(\e^{-1}))^{2(1+p_k)/(3+4p_k)}=o(1),\quad
\end{eqnarray*}
and by \nref{adtom1}, for any $k=0,1,\ldots,K_{\e}$, we have
$$
P_{\e,0}(t_{\e,c_k}>H_\e)\le\frac{E_{\e,0}(\exp(H_\e
t_{\e,c_k}))}{\exp(H_\e^2)}\sim
\exp(H_\e^2/2-H_\e^2)=\exp(-H_\e^2/2)=K_\e^{-2}.
$$
This implies that, for the type I error probability,
$$
\a(\psi_\e)\le \sum_{k=0}^{K_\e} P_{\e,0}(t_{\e,c_k}> H_\e)\le
K_\e^{-1}(1+o(1))\to 0.
$$

Let us evaluate the type II error probability for
$$
\eta\in \Theta_\e(\Sigma)=\bigcup_{p\in\Sigma}\Theta_{\e,p}(r_\e(p)).
$$
There
exists $p$ such that $\eta\in\Theta_{\e,p}(r_\e(p)),\ p_{k}\le p\le
p_{k-1}$.  Observe that
$$
\b_\e(\eta,\psi_\e)\le \min_{0\le k\le
K_\e}P_{\e,\eta}(t_{\e,c_k}\le H_\e).
$$
Denote $h_{\e,c_k}=E_{\e,\eta}(t_{\e,c_k})$. We then have
\begin{equation}\label{adtom2}
h_{\e,c_k}=\e^{-2}\sum_{\nu \in
C_{\nu,k}}w_{\nu,k}\eta_{\nu}^2;\quad
\Var_{\e,\eta}(t_{\e,c_k})=1+4\e^{-2}\sum_{\nu \in
C_{\nu,k}}w_{\nu,k}^2\eta_\nu^2=1+O(h_{\e,c_k}).
\end{equation}
Let us now evaluate
$$
h_{\e,c_k}=\frac{1}{\e^2\sqrt{2F_{p_k}(c_k^2)}}\sum_{\nu \in
C_{\nu,k}}\sigma_\nu^2\eta_{\nu}^2,
$$
where the function
$$ F_p(c)=\sum_{\{\nu:\, a^2_{\nu,p}\le c\}}\sigma_\nu^4
$$
is of the form \nref{Fnew} with asymptotics given by \nref{Fnew} as
well. Observe that $a_{\nu,p_k}=a_{\nu,p}^{p_k/p}$ and, hence,
$$
\sum_{\nu \in C_{\nu,k}}\sigma_\nu^2\eta_\nu^2=\sum_{\nu \in
\Gamma}\sigma_\nu^2\eta_\nu^2-\sum_{\{\nu:\,a_{\nu,p_k} >
c_k\}}\sigma_\nu^2\eta_\nu^2\ge
r_\e^2(p)-c_k^{-2p/p_k}=r_\e^2(p)\left(1-\frac1{r_\e^2(p)c_k^{2p/p_k}}\right).
$$
Because $a(p)^{p\phi(p)/4}$ satisfy the Lipschitz continuity, we
have by \nref{M*1} with $d$ replaced by $D$,
$$
\frac{r_\e(p_k)}{r_\e(p)}=(1+O(\delta_\e))\exp((p\phi(p)-p_k\phi(p_k))(\log(\e^{-1})+(\log(D^{-1}
-\log\log\log(\e^{-1}))/4)),
$$
and
$$
0\le p\phi(p)-p_k\phi(p_k)=\frac34(\phi(p_k)-\phi(p))\le
\frac34\delta_\e\sim\frac{3\log(2)}{4\log(\e^{-1})\log\log(\e^{-1})}.
$$
We get
\begin{equation}\label{adtom3}
 1\le\frac{r_\e(p_k)}{r_\e(p)}\le 2^{3/(4\log\log(\e^{-1}))}(1+o(1))=1+o(1),
\end{equation}
and, moreover, we have
$$
\frac{p}{p_k}=1+\frac{\Delta p}{p_k},\quad 0\le p-p_k=\Delta p\le
p_{k-1}-p_k\asymp \frac1{\log(\e^{-1})\log\log(\e^{-1})}.
$$
Thus, we have
$$
c_k\le c_k^{p/p_k}=c_k\cdot c_k^{\Delta p/p_k}= c_k(1+o(1)),\quad
r_\e(p)c_k^{p/p_k}=2(1+o(1)).
$$
These relations and   \nref{adtom3} imply, for  $\e>0$ small enough,
$$
\sum_{\nu \in C_{\nu,k}}\sigma_\nu^2\eta_\nu^2\ge
r_\e^2(p_k)(1-1/4(1+o(1))=r_\e^2(p_k)(3/4+o(1)).
$$
Therefore, for large enough $D>0$, we have, for all $k=1,2,\ldots,
K_\e$,
\begin{eqnarray}\label{adtom4}
h_{\e,c_k}&\ge& \frac
{r_\e^2(p_k)(3/4+o(1))}{\e^2\sqrt{2F_{p_k}(c_k^2)}}\sim\frac {\
r_\e^2(p_k)}{\e^2c_k^{3/(2p_k)}}\left(\frac{3\sqrt{3}}{8\sqrt{B}}+o(1)\right)\nonumber\\
&\sim&
(Da(p_k)\log\log(\e^{-1}))^{1/2}2^{-3/(2p_k)}\left(\frac{3\sqrt{3}}{8\sqrt{B}}+o(1)\right)\nonumber\\
&>&2H_\e\sim 4\sqrt{\log\log(\e^{-1})}.
\end{eqnarray}

It follows from \nref{adtom2} and \nref{adtom4} that  one has, for
$k$ such that $\eta\in\Theta_{\e,p}(r_\e(p)),\ p_{k}\le p\le
p_{k-1}$,
\begin{eqnarray}\label{adtom6}
P_{\e,\eta}(t_{\e,c_k}\le H_\e)=P_{\e,\eta}(h_{\e,c_k}-t_{\e,c_k}\ge
h_{\e,c_k}-H_\e)\le\frac{\Var_{\e,\eta}(h_{\e,c_k}-t_{\e,c_k})}{(h_{\e,c_k}-H_\e)^2}\nonumber\\
<\frac{4\Var_{\e,\eta}(t_{\e,c_k})}{h_{\e,c_k}^2}=
\frac{4(1+O(h_{\e,c_k}))}{h_{\e,c_k}^2}<
\frac1{4\log\log(\e^{-1})}+o_\eta(1)=o(1).
\end{eqnarray}
Hence, for any $\eta\in \Theta_\e(\Sigma)$, one has
$$
\b_\e(\eta,\psi_\e)<1/(4\log\log(\e^{-1}))+o_\eta(1)=o(1).
$$
$$
\b_\e(\Theta_\e(\Sigma),\psi_\e)=\sup_{\eta\in\Theta_\e(\Sigma)}\b_\e(\eta,\psi_\e)\le
\frac1{4\log\log(\e^{-1})}+\sup_{\eta\in\Theta_\e(\Sigma)}o_\eta(1),
$$
where, by \nref{adtom6}, \nref{adtom5} and \nref{adtom2} for any
$\eta\in\Theta_\e(\Sigma)$
\begin{equation*}
o_\eta(1) = \frac{4\e^{-2}\sum_{\nu \in
C_{\nu,k}}w_{\nu,k}^2\eta_\nu^2}{(\e^{-2}\sum_{\nu \in
C_{\nu,k}}w_{\nu,k}\eta_\nu^2)^2}\le \frac{4\max_{\nu \in
C_{\nu,k}}w_{\nu,k}}{h_{\e,c_k}}<\frac2{H_\e}\sim
(\log\log(\e^{-1}))^{-1/2}.
\end{equation*}
This completes the proof of part (b) of the theorem.

Part (c) of the theorem follows immediately from parts (a) and (b)
of the theorem and \nref{UMFF2}.
The theorem now follows.

\section*{\large Acknowledgements}
We are grateful to the AE and the referee for their useful comments and recommendations.

\label{UMFF}

\end{document}